\def\DATE{\relax}
\magnification=1100
\baselineskip=13pt
\voffset=.75in  
\hoffset=1truein
\hsize=4.5truein
\vsize=7.75truein
\parindent=.166666in
\pretolerance=500 \tolerance=1000 \brokenpenalty=5000

\footline={\hfill{\rm\the\pageno}\hfill\llap{\sevenrm\DATE}}

\def\note#1{%
  \hfuzz=50pt%
  \vadjust{%
    \setbox1=\vtop{%
      \hsize 3cm\parindent=0pt\eightpoints\baselineskip=9pt%
      \rightskip=4mm plus 4mm\raggedright#1%
      }%
    \hbox{\kern-4cm\smash{\box1}\hfil\par}%
    }%
  \hfuzz=0pt
  }
\def\note#1{\relax}

\def\anote#1#2#3{\smash{\kern#1in{\raise#2in\hbox{#3}}}%
  \nointerlineskip}     

\newcount\equanumber
\equanumber=0
\newcount\sectionnumber
\sectionnumber=0
\newcount\subsectionnumber
\subsectionnumber=0
\newcount\snumber  
\snumber=0

\def\section#1{%
  \subsectionnumber=0%
  \snumber=0%
  \equanumber=0%
  \advance\sectionnumber by 1%
  \noindent{\bf \the\sectionnumber .~#1.~}%
}%
\def\subsection#1{%
  \advance\subsectionnumber by 1%
  \snumber=0%
  \equanumber=0%
  \noindent{\bf \the\sectionnumber .\the\subsectionnumber .~#1.~}%
}%
\def\prevs{\the\sectionnumber .\the\subsectionnumber .\the\snumber }

\long\def\Corollary#1{%
  \global\advance\snumber by 1%
  \bigskip
  \noindent{\bf Corollary~\prevs .}%
  \quad{\it#1}%
}%
\long\def\Lemma#1{%
  \global\advance\snumber by 1%
  \bigskip
  \noindent{\bf Lemma~\prevs .}%
  \quad{\it#1}%
}%

\long\def\Proposition#1{%
  \advance\snumber by 1%
  \bigskip
  \noindent{\bf Proposition~\prevs .}%
  \quad{\it#1}%
}%
\long\def\Theorem#1{%
  \advance\snumber by 1%
  \bigskip
  \noindent{\bf Theorem~\prevs .}%
  \quad{\it#1}%
}%
\long\def\Statement#1{%
  \advance\snumber by 1%
  \bigskip
  \noindent{\bf Statement~\prevs .}%
  \quad{\it#1}%
}%
\def\ifundefined#1{\expandafter\ifx\csname#1\endcsname\relax}
\def\labeldef#1{\global\expandafter\edef\csname#1\endcsname{\prevs}}
\def\labelref#1{\expandafter\csname#1\endcsname}
\def\label#1{\ifundefined{#1}\labeldef{#1}\note{$<$#1$>$}\else\labelref{#1}\fi}

\def\preveq{(\the\sectionnumber .\the\subsectionnumber .\the\equanumber)}
\def\neq{\global\advance\equanumber by 1\eqno{\preveq}}

\def\ifundefined#1{\expandafter\ifx\csname#1\endcsname\relax}

\def\equadef#1{\global\advance\equanumber by 1%
  \global\expandafter\edef\csname#1\endcsname{\preveq}%
  \preveq}

\def\equaref#1{\expandafter\csname#1\endcsname}

\font\eightrm=cmr8%
\font\sixrm=cmr6%

\font\eightsl=cmsl8%

\font\eightbf=cmb8%

\font\eighti=cmmi8%
\font\sixi=cmmi6%

\font\eightsy=cmsy8%
\font\sixsy=cmsy6%

\font\eightex=cmex8%
\font\sixex=cmex6%
\font\fiveex=cmex5%

\font\eightit=cmti8%

\font\eighttt=cmtt8%

\font\tenbb=msbm10%
\font\eightbb=msbm8%
\font\sevenbb=msbm7%
\font\sixbb=msbm6%
\font\fivebb=msbm5%
\newfam\bbfam  \textfont\bbfam=\tenbb  \scriptfont\bbfam=\sevenbb  \scriptscriptfont\bbfam=\fivebb%

\font\tencmssi=cmssi10%
\font\sevencmssi=cmssi7%
\font\fivecmssi=cmssi5%
\newfam\ssfam  \textfont\ssfam=\tencmssi  \scriptfont\ssfam=\sevencmssi  \scriptscriptfont\ssfam=\fivecmssi%

\font\tenfrak=cmfrak10%
\font\eightfrak=cmfrak8%
\font\sevenfrak=cmfrak7%
\font\sixfrak=cmfrak6%
\font\fivefrak=cmfrak5%
\newfam\frakfam  \textfont\frakfam=\tenfrak  \scriptfont\frakfam=\sevenfrak  \scriptscriptfont\frakfam=\fivefrak%
\def\frak{\fam\frakfam\tenfrak}%

\font\tenmsam=msam10%
\font\eightmsam=msam8%
\font\sevenmsam=msam7%
\font\sixmsam=msam6%
\font\fivemsam=msam5%

\def\bb{\fam\bbfam\tenbb}%

\def\hexdigit#1{\ifnum#1<10 \number#1\else%
  \ifnum#1=10 A\else\ifnum#1=11 B\else\ifnum#1=12 C\else%
  \ifnum#1=13 D\else\ifnum#1=14 E\else\ifnum#1=15 F\fi%
  \fi\fi\fi\fi\fi\fi}
\newfam\msamfam  \textfont\msamfam=\tenmsam  \scriptfont\msamfam=\sevenmsam  \scriptscriptfont\msamfam=\fivemsam%
\def\msam{\msamfam\tenmsam}%
\mathchardef\leq"3\hexdigit\msamfam 36%
\mathchardef\geq"3\hexdigit\msamfam 3E%

\newdimen\bigsize  \bigsize=8.5pt
\newdimen\Bigsize  \Bigsize=11.5pt
\newdimen\biggsize \biggsize=14.5pt
\newdimen\Biggsize \Biggsize=17.5pt

\catcode`\@=11
\def\big#1{{\hbox{$\left#1\vbox to\bigsize{}\right.\n@space$}}}
\def\Big#1{{\hbox{$\left#1\vbox to\Bigsize{}\right.\n@space$}}}
\def\bigg#1{{\hbox{$\left#1\vbox to\biggsize{}\right.\n@space$}}}
\def\Bigg#1{{\hbox{$\left#1\vbox to\Biggsize{}\right.\n@space$}}}
\catcode`\@=12

\def\eightpoints{%

\def\rm{\fam0\eightrm}%
\textfont0=\eightrm   \scriptfont0=\sixrm   \scriptscriptfont0=\fiverm%
\textfont1=\eighti    \scriptfont1=\sixi    \scriptscriptfont1=\fivei%
\textfont2=\eightsy   \scriptfont2=\sixsy   \scriptscriptfont2=\fivesy%
\textfont3=\eightex   \scriptfont3=\sixex   \scriptscriptfont3=\fiveex%
\textfont\itfam=\eightit  \def\it{\fam\itfam\eightit}%
\textfont\slfam=\eightsl  \def\sl{\fam\slfam\eightsl}%
\textfont\ttfam=\eighttt  \def\tt{\fam\ttfam\eighttt}%
\textfont\bffam=\eightbf  \def\bf{\fam\bffam\eightbf}%

\textfont\frakfam=\eightfrak  \scriptfont\frakfam=\sixfrak \scriptscriptfont\frakfam=\fivefrak  \def\frak{\fam\frakfam\eightfrak}%
\textfont\bbfam=\eightbb      \scriptfont\bbfam=\sixbb     \scriptscriptfont\bbfam=\fivebb      \def\bb{\fam\bbfam\eightbb}%
\textfont\msamfam=\eightmsam  \scriptfont\msamfam=\sixmsam \scriptscriptfont\msamfam=\fivemsam  \def\msam{\msamfam\eightmsam}

\bigsize=6.8pt
\Bigsize=9.2pt
\biggsize=11.6pt
\Biggsize=14pt

\rm%
}

\mathchardef\lsim"3\hexdigit\msamfam 2E%
\mathchardef\gsim"3\hexdigit\msamfam 26%

\def\ds{\displaystyle}
\long\def\DoNotPrint#1{\relax}

\def\qed{{\vrule height .9ex width .8ex depth -.1ex}}

\def\phi{\varphi}

\def\calF{{\cal F}}

\def\NN{{\bb N}\kern .5pt}

\pageno=1

\centerline{\bf AN ELEMENTARY APPROACH}
\centerline{\bf TO EXTREME VALUES THEORY}
\bigskip
 
\centerline{Ph.\ Barbe}
\centerline{CNRS, France}
 
{\narrower
\baselineskip=9pt\parindent=0pt\eightpoints

\bigskip

{\bf Abstract.} This note presents a rather intuitive approach to
extreme value theory. This approach was devised mostly for pedagogical
reason.

\bigskip


\bigskip
 

}

\bigskip

\noindent {\bf 1. Introduction.} 
The purpose of this note is to present a rather elementary approach to
some results in extreme values theory. The main proof was designed mostly 
for pedagogical reasons so that it could be taught at a very intuitive 
level. In particular, the main
result does not use regular variation theory or the concept of type of
a distribution.

To recall what extreme value theory is about, let $M_n$ be the maximum 
of $n$ independent real
random variables all having the same distribution function $F$.
Extreme value theory grew from the search (now completed) for necessary and 
sufficient condition for $M_n$ linearly normalized to have a nondegenerate 
limiting distribution.
In particular, one says that $F$ belongs to a domain of max-attraction if 
there exist deterministic 
sequences $(a_n)_{n\geq 1}$ and $(b_n)_{n\geq 1}$ such that the distribution
of $(M_n-b_n)/a_n$ converges to a nondegenerate limit. This makes the root of
what we
call here {\it linear} extreme value theory. In contrast, {\it nonlinear}
extreme value theory seeks sequences of deterministic and possibly nonlinear
functions $(g_n)_{n\geq 1}$ such that the distribution of $g_n(M_n)$ converges
to a nondegenerate limit as $n$ tends to infinity. Note that in this
context, it is rather natural to restrict each $g_n$ to be monotone, and,
then without any loss of generality, to be nondecreasing, eventually by
replacing $g_n$ by $-g_n$.

Of essential importance for both the linear and nonlinear extreme value
theory, the quantile function pertaining to the distribution 
function $F$ is defined as
$$
  F^\leftarrow(u)=\inf\{\, x\,:\, F(x)>u\,\} \, .
$$
It is c\`adl\`ag, that is right continuous with left limits, as well as 
nondecreasing.

Note that in linear extreme value theory, only three possible limiting
distribution can arise. In contrast, any nondegenerate limiting distribution 
can arise using nonlinear normalization. This can be seen very easily as 
follows. Consider an arbitrary distribution function $G$ and consider
$F$ to be the uniform distribution over $(0,1)$. 
Define $g_n(x)=G^\leftarrow(e^{ -n(1-x)})$. A direct 
calculation shows that the distribution of $n(1-M_n)$ converges to the 
standard exponential one. Consequently, $g_n(M_n)$ has limiting
distribution $G$.

\bigskip

\noindent{\bf 2. Nonlinear extreme value theory.} The following result
characterizes all possible monotone transformations $g_n$ of $M_n$ such 
that $g_n(M_n)$ has a nondegenerate limiting distribution.

\bigskip

\noindent{\bf Theorem 2.1.}{\it\quad
  Let $(g_n)_{n\geq 1}$ be a sequence of nondecreasing
  functions on the real line.
  The following are equivalent:

  \smallskip

  \noindent (i) The distribution of $g_n(M_n)$ converges 
  to a nondegenerate limit.

  \smallskip

  \noindent (ii) The sequence of functions 
  $x\in (0,\infty) \mapsto g_n\circ F^\leftarrow (1-x/n)$, $n\geq 1$, 
  converges almost everywhere to a nonconstant limit.

  \smallskip

  \noindent In this case, writing $\omega$ for a standard exponential 
  random variable,
  and $h$ for the limiting function involved in (ii), the limiting 
  distribution function of $g_n(M_n)$ is that of $h(\omega)$. Moreover,
  $h$ is continuous almost everywhere.
}

\bigskip

It is easy to see from its proof that Theorem 2.1 still holds if one 
replaces the full sequence $(n)_{n\geq 1}$ by a subsequence $(n_k)_{k\geq 1}$.
Thus, the same result, considering now subsequences, applies 
to so-called partial domain of attraction.
The technique used in the proof also shows that assertion (ii) in Theorem 2.1
is equivalent to
$$
  \lim_{\epsilon\to 0} g_{\lfloor 1/\epsilon\rfloor}\circ F^\leftarrow 
  (1-\epsilon x)
  \eqno{(2.1)}
$$
exists almost everywhere and is nonconstant on $(0,\infty)$.

\bigskip

Because of the pedagogical motivation of this note, we give a complete
proof of Theorem 2.1 as far as the probabilistic arguments are
concerned. We will need some known auxiliary results which we state as 
lemmas and whose proofs are given for pedagogical reasons but deferred 
to an appendix.

Our first lemma is the so-called
quantile transform which consists of the following known result.

\bigskip

\noindent{\bf Lemma 2.2.}\quad {\it 
  Let $U$ be a random variable having a uniform distribution over $(0,1)$.
  The random variable $F^\leftarrow (U)$ has distribution function $F$.
  }

\bigskip

The second lemma collects two elementary facts on convergence of sequences
of functions, the first assertion being not much more than a restatement
of Helly's theorem (see Feller, 1970, \S VIII.6), and the whole lemma
being exercise 13 in chapter 7 of Rudin's (1986) {\it Principles of 
Mathematical Analysis}.

\bigskip

\noindent{\bf Lemma 2.3.}\quad{\it
  (i) A uniformly locally bounded sequence of nonincreasing 
  functions has an almost everywhere convergent subsequence whose limit
  is continuous almost everywhere.

  \smallskip

  \noindent (ii) A family of nonincreasing functions
  which converges almost everywhere to a continuous limit converges 
  everywhere and locally
  uniformly.
}

\bigskip

\noindent{\bf Proof of Theorem 2.1.} Let $U_n$
be the maximum of $n$ independent random variables uniformly distributed
on $(0,1)$. A direct calculation shows that
$$
  P\{\, U_n\leq x\,\} = x^n \, .
  \eqno{(2.2)}
$$
Let $\omega$ be a random variable
having the standard exponential distribution. Note that $e^{-\omega}$
is uniformly distributed over $(0,1)$. Thus, (2.2) implies that $U_n$
has the same distribution as $e^{-\omega/n}$. 
Using the quantile transform, that
is Lemma 2.2, we see that the distribution of $M_n$ is that of 
$F^\leftarrow(U_n)$, that is, that
of $F^\leftarrow(e^{-\omega/n})$.  Therefore, for $g_n(M_n)$ to
have a nondegenerate limiting distribution,
it is necessary and sufficient that the distribution of
$$
  h_n(\omega)=g_n\circ F^\leftarrow(e^{-\omega/n})
  \eqno{(2.3)}
$$
converges as $n$ tends to infinity. The intuition behind our proof is that
if this convergence holds then it holds almost surely because the random
variable $\omega$ does not depend on $n$. Thus, we will first consider the
assertion%
\setbox1=\vbox{\narrower\narrower\narrower\noindent
  the sequence $(h_n)_{n\geq 1}$ converges almost everywhere to a
  limit which is nonconstant on $(0,\infty)$.
  }\dp1=0pt\ht1=0pt%
\vskip-20pt
$$
  \box1
  \eqno{(2.4)}
$$
\vskip 15pt

\noindent{\it Proof that (2.4) implies (i)}. If (2.4) holds,
call $h$ the limit of the sequence $(h_n)_{n\geq 1}$. 
Since $\lim_{n\to\infty} h_n=h$ almost everywhere,
the distribution of the random variable 
$h_n(\omega)$ converges to that of $h(\omega)$ as $n$ tends to infinity.
Since $h$ is nonincreasing and is not constant, there exists a real number
$a$ such that $h(0,a)\cap h(a,\infty)$ is empty. This implies that the
random variable $h(\omega)$ is nondegenerate.

\noindent{\it Proof that (i) implies (2.4)}. Let $G$ be the nondegenerate
limiting distribution involved in (i).
In order to prove that the sequence $h_n$ defined in (2.3) converges, 
we first show that it
satisfies the assumptions of Lemma 2.3. Note that each function $h_n$ is
nonincreasing.

\bigskip

\noindent{\bf Lemma 2.4.} \quad{\it
  The sequence $(h_n)_{n\geq 1}$ is locally uniformly bounded on $(0,\infty)$.
  }

\bigskip

\noindent{\bf Proof.} Let $[\,a,b\,]$ be a bounded interval in $(0,\infty)$. 
Seeking a contradiction,
assume that the sequence $(h_n)_{n\geq 1}$ is not bounded on $[\, a,b\,]$.
Then, we can extract a subsequence $(\omega_k)_{k\geq 1}$ in $[\,a,b\,]$ and
a subsequence $(n_k)_{k\geq 1}$ such that $h_{n_k}(\omega_k)$ tends to either
$+\infty$ or $-\infty$. Assume first 
that $\lim_{k\to\infty} h_{n_k}(\omega_k)=+\infty$. Since $h_{n_k}$ is 
nonincreasing, $\lim_{k\to\infty}h_{n_k}(a)=+\infty$. Therefore, for any $M$
positive and any $k$ large enough,
$$\eqalign{
  1-e^{-a}
  =P\{\, \omega \leq a\,\}
  &{}\leq  P\{\, h_{n_k}(\omega)\geq h_{n_k}(a)\,\}\cr
  &{}\leq P\{\, h_{n_k}(\omega)\geq M\,\} \, .\cr
 }
$$
Taking limit as $k$ tends to infinity we obtain $1-e^{-a}\leq 1-G(M-)$. 
Since $M$ is arbitrary large, this yields $1-e^{-a}
\leq 0$, which is the desired contradiction.

If we assume that $\lim_{k\to\infty}h_{n_k}(\omega_k)=-\infty$, then
$\lim_{k\to\infty}h_{n_k}(b)=-\infty$. Therefore, for
any $M$ negative and any $k$ large enough,
$$\eqalign{
  e^{-b}
  = P\{\, \omega>b\,\}
  &{}\leq P\{\, h_{n_k}(\omega)\leq h_{n_k}(b)\,\} \cr
  &{}\leq P\{\, h_{n_k}(\omega)\leq M \,\} \, .\cr
  }
$$
Taking limit as $k$ tends to infinity yields $e^{-b}\leq G(M)$, and since $M$
is arbitrary, $e^{-b}\leq 0$, which is a contradiction.\hfill\qed

\bigskip

From Lemmas 2.3 and 2.4 we deduce that we can find a subsequence $h_{n_k}$
which converges almost everywhere to a limit $h$, and, moreover, this 
limit is nonincreasing. But then,
$$
  P\{\, h(\omega)\leq x\,\} = G(x) \, .
$$
It follows that $h$ is unique almost everywhere and that any convergent
subsequence of $(h_n)_{n\geq 1}$ converges to $h$. Then, Lemma 2.4 implies that
the sequence $(h_n)_{n\geq 1}$ converges almost everywhere to $h$.

\noindent{\it Equivalence between (2.4) and Theorem 2.1.ii}.
We consider the sequence of functions
$$ 
  \tilde h_n(\omega)=g_n\circ F^\leftarrow (1-\omega/n) \, .
$$
Since $e^{-\omega/n}\geq 1-\omega/n$, we see that $h_n\geq\tilde h_n$. 
For any fixed $\omega$ and any $n$ large enough, $e^{-\omega/n(1-\epsilon)}\leq
1-\omega/n$. Therefore, for $n$ large enough,
$\tilde h_n(\omega)\geq h_n\bigl( \omega/(1-\epsilon)\bigr)$.
If (ii) holds the above inequalities comparing $h_n$ and $\tilde h_n$ show
that
$$
  h\bigl( \omega/(1-\epsilon)\bigr)
  \leq \liminf_{n\to\infty} \tilde h_n(\omega)
  \leq \limsup_{n\to\infty} \tilde h_n(\omega)
  \leq h(\omega)
$$
almost everywhere. If $\omega$ is a continuity point of $h$, then
$h\bigl( \omega/(1-\epsilon)\bigr)$ tends to $h(\omega)$ as $\epsilon$
tends to $0$,  
and, consequently, $\tilde h_n(\omega)$ converges to $h(\omega)$.

Conversely, if (2.4) holds, the limiting function $h$ is monotone and locally
bounded. Hence it has at most countable many discontinuities and it is almost
everywhere continuous. The same 
bound relating $h_n$ and $\tilde h_n$ show that $h_n$ converges almost
everywhere to $h$, which is (ii).

\noindent{\it Equivalence between (2.1) and (ii)}.
Clearly, if (2.1) holds then assertion (ii) of Theorem 2.1 holds. 
To prove the converse
implication, let $x$ be a point of continuity of $h$ such that
$$
  \lim_{n\to\infty} g_n\circ F^\leftarrow (1-x/n) = h(x) \, .
$$
Let $n$ be the integer part of $1/\epsilon$, so that $1/(n+1)<
\epsilon\leq 1/n$. For any fixed $\eta$, provided that $\epsilon$ is small
enough,
$$
  F^\leftarrow (1-x/n) 
  \leq F^\leftarrow(1-\epsilon x)
  \leq F^\leftarrow \bigl(1-(x-\eta)/n\bigr) \, .
$$
In particular,
$$
  g_n\circ F^\leftarrow (1-x/n)
  \leq g_{\lfloor 1/\epsilon\rfloor}\circ F^\leftarrow (1-\epsilon x)
  \leq g_n \circ F^\leftarrow \bigl( 1-(x-\eta)/n\bigr) 
  \, .
$$
Taking limit as $\epsilon$ tends to $0$ and then limit as $\eta$ tends to $0$
and using that $x$ is a continuity point of $h$,
$$\displaylines{\qquad
  h(x)
  \leq \liminf_{\epsilon\to 0} g_{\lfloor 1/\epsilon\rfloor}
   \circ F^\leftarrow (1-\epsilon x)
  \hfill\cr\hfill
  \leq \limsup_{\epsilon\to 0}g_{\lfloor 1/\epsilon\rfloor}
   \circ F^\leftarrow(1-\epsilon x)
   \leq h(x) \, . \qquad\cr}
$$
This proves (2.1).\hfill\qed

\bigskip

\noindent{\bf 3. Application to linear extreme value theory.}
The purpose of this section is to show how some classical
results can be derived from Theorem 2.1. We mostly restrict ourself to
the following result, due to de Haan (1970), which characterizes the 
belonging to a domains of attraction.

\bigskip

\noindent{\bf Theorem 3.1}~(de Haan, 1970).\quad {\it 
  A distribution function $F$ belongs to a domain of max-attraction if and
  only if for any
  $$
    \lim_{\epsilon\to 0} 
    { F^\leftarrow (1-\epsilon u)-F^\leftarrow(1-\epsilon)
      \over
      F^\leftarrow (1-\epsilon v)-F^\leftarrow(1-\epsilon) }
    \quad \hbox{ exists}
    \eqno{(3.1)}
  $$
  for almost all $u$ and $v$.
}

\bigskip

\noindent{\bf Remark.} Theorem 3.1 does not state the classical convergence of
type result, namely that there are only three
possible types of limiting distribution. This can be recovered by the following
known argument. For any real number $\rho$, define the function
$$
  k_\rho(u)=\cases{ {\ds u^\rho-1\over\ds \rho} & if $\rho\not= 0$,\cr
                    \noalign{\vskip 3pt}
                    \log u               & if $\rho=0$.\cr}
$$
It can be shown (see Bingham, Goldie and Teugels,
1989, chapter 3, or the appendix to this paper which reproduces their argument
with an extra monotonicity assumption which holds here and leads to substantial
simplifications)
that the limit in (3.1) is
necessarily of the form $k_\rho(u)/k_\rho(v)$ for some real number $\rho$. 
Then, taking 
$$
  a_n=F^\leftarrow(1-2/n)-F^\leftarrow(1-1/n)
  \hbox{ and }
  b_n=F^\leftarrow(1-1/n) \, ,
  \eqno{(3.2)}
$$
we obtain that the distribution of $(M_n-b_n)/a_n$ converges to that
of $k_\rho(\omega)/k_\rho(2)$. An explicit calculation of the limiting
distribution is then easy, and the discussion according to the position of
$\rho$ with respect to $0$ (larger, smaller or equal) yields the classical
three types.

\bigskip

\noindent{\bf Proof of Theorem 3.1.} We mostly present the part of the proof 
related to Theorem 2.1.

\noindent{\it Necessity.} Assume that $F$ belongs to a domain of attraction.
Consider the norming constants $(a_n)_{n\geq 1}$ and $(b_n)_{n\geq 1}$, as
well as the functions $g_n(u)=(u-b_n)/a_n$. 
Define $h_{1/\epsilon}(x)=g_{\lfloor 1/\epsilon\rfloor}\circ 
F^\leftarrow (1-\epsilon x)$.
Theorem 2.1.ii in its formulation (2.1) asserts that $h_{1/\epsilon}$ converges
almost everywhere to some function $h$ as $\epsilon$ tends to $0$. It 
follows that for almost $u,v,x,y$ for which $h(v)$ and $h(y)$ are distinct,
$$\eqalignno{
  {h(u)-h(x)\over h(v)-h(y)}
  &{}= \lim_{\epsilon\to 0} { h_{1/\epsilon}(u)-h_{1/\epsilon}(x)\over 
                              h_{1/\epsilon}(v)-h_{1/\epsilon}(y) } \cr
  &{}=\lim_{\epsilon\to 0} 
   { F^\leftarrow(1-\epsilon u)-F^\leftarrow(1-\epsilon x) \over
     F^\leftarrow(1-\epsilon v)-F^\leftarrow(1-\epsilon y) } \, .
  &(3.3)\cr
  }
$$
This is not quite (3.1) since, a priori, we may not be able to choose $x$ 
and $y$ to 
be $1$. An extra regular variation theoretic argument, essentially explained 
in Bingham, Goldie and Teugels (1989, chapter 3) is then needed. For the sake
of completeness and given the pedagogical nature of this note, we develop
this argument in the appendix.

\noindent{\it Sufficiency.} If (3.1) holds then it holds everywhere and 
locally uniformly and the limit is of the form $k_\rho(u)/k_\rho(v)$
---~see Bingham,
Goldie and Teugels, 1989, Chapter 3; or, alternatively, use the regular
variation theoretic argument in the appendix.
Taking $a_n$ and $b_n$ as in (3.2),
this implies that $\bigl(F^\leftarrow (e^{-\omega/n})-b_n\bigr)/a_n$ has
a limit $k_\rho(\omega)/k_\rho(2)$ as $n$ tends to infinity.
This implies (see the representation for $M_n$ in the proof of Theorem 2.1), 
that the distribution of $(M_n-b_n)/a_n$ converges to a nondegenerate 
limit.\hfill\qed

\bigskip

\noindent{\bf 4. On the maximum of geometric random variables.} In this section
we consider the maximum $M_n$ of $n$ independent random variables all having
a geometric distribution. With the notation of section 1 and writing
$\lfloor\cdot\rfloor$ for the integer part, the underlying
distribution function is
$$
  F(t)=(1-p)\sum_{0\leq i\leq t} p^i = 1-p^{\lfloor t+1\rfloor}
$$
for some  $p$ between $0$ and $1$.
It is known (see e.g.\ Resnick, 1987, \S 1.1, example following Corollary
1.6) that there are no sequences $(a_n)_{n\geq 1}$ and $(b_n)_{n\geq 1}$
such that the distribution of $(M_n-b_n)/a_n$ has a nondegenerate limiting
distribution. In other words, it is not possible to find linear normalizations
or a sequence of deterministic affine functions $(g_n)_{n\geq 1}$ such that
the distribution of $g_n(M_n)$ converges to a nondegenerate limit. A natural
question is then: can we find a sequence of nonlinear functions 
$(g_n)_{n\geq 1}$ such that the distribution of $g_n(M_n)$ has a nondegenerate
limit? The next proposition shows that under the additional requirement that 
each $g_n$ is monotone, the answer is negative. Hence, in some sense, there 
is no good alternative to using subsequences and partial domain of attraction
---~see also the remark following the proof. The same result can
be obtained in combining theorems 1.5.1 and 1.7.13 in Leadbetter, Lindgren
and Rootz\'en (1983). 

\bigskip

\noindent{\bf Proposition 4.1.}\quad{\it
  There is no deterministic sequence of nondecreasing functions 
  $(g_n)_{n\geq 1}$ such that the distribution of $g_n(M_n)$ has a 
  nondegenerate limit.
  }

\bigskip

\noindent{\bf Proof.} The proof is by contradiction and relies on Theorem 2.1.
It also uses the following facts, stated as a lemma, which is a 
classical exercise in 
analytic number theory (see Hlawka, Schoi\ss engeier, Taschner, 1986,
Chapter 2, exercise 8) and whose proof is in the appendix. We write
$\calF(\cdot )$ for the fractional part, that is 
$\calF(x)=x-\lfloor x\rfloor$.

\bigskip

\noindent{\bf Lemma 4.2.}\quad{\it
  For any positive real number $\theta$, the sequence\hfill\break
  $\bigl(\calF(\theta\log n)\bigr)_{n\geq 1}$ is dense in $[\,0,1\,]$.
}

\bigskip

In order to prove Proposition 3, and seeking a contradiction, assume that
there exists a deterministic sequence $(g_n)_{n\geq 1}$ of nondecreasing
functions such that the distribution of $g_n(M_n)$ has a nondegenerate limit.
Theorem 2.1 implies that $g_n\circ F^\leftarrow(1-x/n)$ has a limit almost
everywhere, $h(x)$, which is nonconstant and nonincreasing.

We first calculate the quantile function
$$\eqalign{
  F^\leftarrow(1-u)
  &{}=\inf\{\, t\,:\, p^{\lfloor t+1\rfloor}<u\,\} \cr
  &{}=\inf\Bigl\{\, t\,:\, \lfloor t+1\rfloor >{\log u×\over\log p}\,\Bigr\} \cr
  &{}=\Bigl\lfloor{\log u\over\log p}\Bigr\rfloor \, .\cr
  }
$$
In particular,
$$
  g_n\circ F^\leftarrow (1-x/n)
  =g_n\Bigl( \Bigl\lfloor -{\log n\over\log p} 
  + {\log x\over\log p}\Bigr\rfloor
  \Bigr) \, .
$$
Set $\theta=-1/\log p$ and $y=\log x/\log p$. We then have
$$
  \lim_{n\to\infty} g_n(\lfloor\theta\log n+y\rfloor)=h(p^y) \, .
  \eqno{(4.1)}
$$
Define the functions
$$
  k_n(u)=g_n(u+\lfloor\theta\log n\rfloor)
  \quad \hbox{ and } \quad
  k(y)=h(p^y) \, .
$$
Equality (4.1) is equivalent to
$$
  \lim_{n\to\infty} 
  k_n(\lfloor \theta\log n+y\rfloor-\lfloor\theta\log n\rfloor)
  = k(y) \, .
$$
The advantage of this equality compared to (4.1) is that for fixed $y$ the
argument of $k_n$ remains of order $1$, while the argument of $g_n$ in (4.1)
tends to infinity with $n$. Clearly, the argument of $k_n$, that is,
$\lfloor \theta\log n+y\rfloor-\lfloor\theta\log n\rfloor$, is an integer.
It is equal to an integer $q$ if and only if $\lfloor\theta\log n+y\rfloor
=q+\lfloor\theta\log n\rfloor$, that is, if
$$
  q+\lfloor\theta\log n\rfloor
  \leq \theta\log n +y
  < q+\lfloor \theta\log n\rfloor  +1\, ,
$$
or, equivalently,
$$
  q\in \calF(\theta\log n)+(y-1,y\,] \, .
$$
Moreover, if this inequality holds then
$$
  k_n(\lfloor\theta\log n+y\rfloor-\lfloor\theta\log n\rfloor)
  = k_n(q)\, ,
$$ 
and therefore $\lim_{n\to\infty}k_n(q)=h(y)$. 

Since $h$ is nonconstant and is noincreasing, we can find $y_1$ and
$y_2$ such that $y_1<y_2<y_1+1$ and $h(y_2)<h(y_1)$.  Note that for any
integer $n$, the intervals $\calF(\theta\log n)+(y_1-1,y_1\,]$ and 
$\calF(\theta\log n)+(y_2-1,y_2\,]$ have a nonempty intersection equal to the
interval $\calF(\theta\log n)+(y_2-1,y_1\,]$. Let $\epsilon$ be
a positive real number such that $2\epsilon<h(y_1)-h(y_2)$. Since the
sequence $\bigl(\calF(\theta\log n)\bigr)_{n\geq 1}$ is dense in
$[\,0,1\,]$, there exists infinitely many $n$ such that the intersections
$\calF(\theta\log n)+(y_2-1,y_1\,]$ contain the same integer $q$. For 
those $n$ sufficiently large, we then have
$$
  |k_n(q)-h(y_1)|<\epsilon
  \quad\hbox{ and }\quad
  |k_n(q)-h(y_2)|<\epsilon \, ,
$$
which forces $|h(y_1)-h(y_2)|<2\epsilon$ and contradicts our choice of
$\epsilon$.\hfill\qed

\bigskip

\noindent {\bf Remark.} The proof shows in fact a little more, namely, that if
there exists a deterministic sequence $(g_n)_{n\geq 1}$ of nondecreasing
functions and if there exists a subsequence $n_k$ such that the
distribution of $g_{n_k}(M_{n_k})$ converges to a nondegenerate limit as $k$
tends to infinity, then it is necessary that the sequence 
$\bigl(\calF(\log n_k)\bigr)$ is not dense in $[\,0,1\,]$. This forces the 
sequence $(n_k)_{k\geq 1}$ to avoid a set of the form
$\cup_{q\in\NN}(e^q[\,e^x,e^y\,])$ for some $0<x<y$, and hence forces that
sequence to contain gaps which grow at least geometrically. 

\bigskip

\noindent{\bf Appendix.}

\medskip

{

\noindent {\bf Proof of Lemma 2.1.} If $s>F^\leftarrow (U)$ then $F(s)>U$. 
Therefore,
$$
  P\{\, F^\leftarrow (U)<s\,\}
  \leq P\{\, U\leq F(s)\,\}
  = F(s) \, .
$$
Since distribution functions are right continuous, this implies
$$
  P\{\, F^\leftarrow (U)\leq s\,\} \leq F(s) \, .
$$
If $s<F^\leftarrow (U)$ then $F(s)\leq U$. Therefore,
$$\eqalignno{
  P\{\, F^\leftarrow (U)\leq s\,\}
  & {}= 1-P\{\, F^\leftarrow (U)\geq s\,\} \cr
  & {}\geq 1-P\{\, U\geq F(s)\,\} \cr
  & {}=F(s) \, . &
  \qed\cr}
$$

\bigskip

\noindent{\bf Proof of Lemma 2.2.} (i) A quick proof consists in considering 
that up to replacing nonincreasing by nondecreasing
such sequence
defines a sequence of measure on the compact sets $[\,0,\infty\,]$ as
well as $[\,-\infty,0\,]$ and use
Prohorov's theorem (see Billingsley, 1968, Theorem 6.1). A more pedestrian
approach is to spell out the arguments as follows. Let $(f_n)_{n\geq 1}$
be a sequence as in the lemma and let $(x_k)_{k\geq 1}$ be a sequence of
numbers dense in the real line. Since the sequence 
$\bigl(f_n(x_1)\bigr)_{n\geq 1}$ is bounded, we can find an increasing function
$\phi_1$ mapping $\NN$ into itself such that 
$\bigl(f_{\phi_1(n)}(x_1)\bigr)_{n\geq 1}$ converges. Suppose that we
constructed an increasing function $\phi_k$ from $\NN$ into itself. We 
construct $\phi_{k+1}$ by requiring that it is increasing, maps $\NN$ into
$\phi_k(\NN)$, that is, $\bigl(\phi_{k+1}(n)\bigr)_{n\geq 1}$ is a 
subsequence of $\bigl(\phi_k(n)\bigr)_{n\geq 1}$, and $\bigl( f_{\phi_{k+1}(n)}
(x_{k+1})\bigr)_{n\geq 1}$ converges. Then, for any fixed $k$ the sequence
$\bigl( f_{\phi_n(n)}(x_k)\bigr)_{n\geq 1}$ converges, and we write 
$\underline f(x_k)$ its limit. Since the functions $f_n$ are nonincreasing,
so are the function $f_{\phi(n)}(n)$ and so is $\underline f$ on the 
set $(x_k)_{k\geq 1}$. Moreover, $\underline f$ is locally bounded. We extend
$\underline f$ to a function $f$ defined on the whole real line by setting
$$
  f(x)
  = \lim_{x_k\downarrow x}\underline f(x_k) 
  = \sup\{\, f(x_k)\,:\, x_k\geq x\,\} \, .
$$
Since $\underline f$ is nonincreasing on the set $(x_k)_{k\geq 1}$, the 
function $f$ is nonincreasing on the real line. Consequently, it has left limit
everywhere. It is right continuous because if $x<y<x+\epsilon$ then we can
find $x_k$ and $x_\ell$ such that $x<x_k<y<x_\ell<x+\epsilon$, which 
implies $f(x_k)\geq f(y)\geq f(x_\ell)$; hence
$$
  f(x)
  =\lim_{x_k\downarrow x}\underline f(x_k)
  \geq \limsup_{y\downarrow x}f(y)
  \geq \liminf_{y\downarrow x} f(y)
  \geq \lim_{x_\ell\downarrow x} \underline f(x_\ell)
  = f(x) \, .
$$
This proves that $f$ is c\`adl\`ag. Since it is locally bounded, it has
countable many discontinuity points. Hence, almost every real number is
a continuity point of $f$. Let $x$ be a continuity point of $f$
and let us prove that $\bigl( f_{\phi_n(n)}(x)\bigr)_{n\geq 1}$ converges to
$f(x)$. Indeed, if $x_i<x<x_k$, then $f_{\phi_n(n)}(x_i)\geq f_{\phi_n(n)}(x)
\geq f_{\phi_n(n)}(x_k)$. Thus, taking limits as $n$ tends to infinity,
$$
  f(x_i)
  \geq \limsup_{n\to\infty} f_{\phi_n(n)}(x)
  \geq \liminf_{n\to\infty} f_{\phi_n(x)}(x)
  \geq f(x_k) \, .
$$
Since $x$ is a continuity point of $f$, taking the limits as $x_i$ and $x_k$ 
tend to $x$ shows that $\lim_{n\to\infty}f_{\phi_n(n)}(x)=f(x)$. This implies
that the subsequence $\bigl( f_{\phi_n(n)} \bigr)_{n\geq 1}$ converges
almost everywhere.

\noindent (ii) Consider an interval $[\,a,b\,]$. 
Let $\eta$ be a positive real number.
The function $\Delta$ being continuous, it is uniformly continuous on 
$[\,a,b\,]$. Moreover, since $\Delta_\epsilon$ is nonincreasing, so is
$\Delta$. Thus, we can find points $a=a_0<a_1<\ldots <a_k=b$ such that
for all $0\leq i\leq k$,
$$
  0\leq \Delta(a_i)-\Delta(a_{i+1})\leq \eta
  \qquad\hbox{and}\qquad
  \lim_{\epsilon\to 0} \Delta_\epsilon(a_i)=\Delta(a_i) \, .
$$
Provided $\epsilon$ is small enough, $|\Delta_\epsilon (a_i)-\Delta(a_i)|
\leq \eta$ for any $1\leq i\leq k$. Consequently, if $x$ is between $a_k$
and $a_{k+1}$,
$$\eqalignno{
  |\Delta_\epsilon(x)-\Delta(x)|
  &{}\leq |\Delta_\epsilon(x)-\Delta_\epsilon(a_{k+1})| 
   + |\Delta_\epsilon(a_{k+1})-\Delta(a_{k+1})| \cr
  &\hskip 3cm {}+ \Delta(a_{k+1})-\Delta(x) \cr
  &{}\leq \Delta_\epsilon(a_k)-\Delta_\epsilon (a_{k+1})+2\eta \cr
  &{}\leq |\Delta_\epsilon(a_k)-\Delta(a_k)| 
   +\Delta_\epsilon(a_{k+1})-\Delta(a_{k+1})| \cr
  &\hskip 3cm {}+ \Delta(a_k)-\Delta(a_{k+1})
   + 2\eta \cr
  &{}\leq 5 \eta \, .
  &\qed\cr
  }
$$

\bigskip

\noindent{\bf Proof of Lemma 4.2.} Consider an interval $[\, x,y\,]$ in 
$[\,0,1\,]$. Let $q$ be an integer. If $q+x\leq\theta\log n\leq q+y$
then $\calF(\theta\log n)$ belongs to $[\, x,y\,]$. Such $n$ exists if
the interval $[\, e^{(q+x)/\theta},e^{(q+y)/\theta}\,]$ contains an
integer. The length of this interval is
$e^{(q+x)/\theta}(e^{(y-x)/\theta}-1)$ and tends to infinity with
$q$. Hence, this interval contains an integer whenever $q$ is large 
enough.\hfill\qed

\bigskip

\noindent{\bf Regular variation theoretic argument for the proof of Theorem 3.1.}
We first give the extra argument needed to prove the necessity part of Theorem
3.1.

Let $v$ and $y$ be such that $h(v)$ and $h(y)$ are distinct, $h_{1/\epsilon}
(v)$ and $h_{1/\epsilon}(y)$ converge to $h(v)$ and $h(y)$
respectively as $\epsilon$ tends to $0$. Define the function
$$
  r(1/\epsilon)=F^\leftarrow(1-\epsilon v)-F^\leftarrow(1-\epsilon y)\,.
$$
Writing $\tilde h(u)=\bigl(h(u)-h(x)\bigr)/\bigl( h(v)-h(y)\bigr)$, (3.3) 
asserts that
$$
  \lim_{\epsilon\to 0} 
  { F^\leftarrow (1-\epsilon u)-F^\leftarrow(1-\epsilon y)
    \over r(1/\epsilon) } =\tilde h(u)
$$
for almost all $u$. In particular, for almost all $u$ and $w$,
$$\displaylines{\qquad
  \tilde h(uw)=\lim_{\epsilon\to 0} 
  { F^\leftarrow (1-\epsilon uw)-F^\leftarrow(1-\epsilon wy)\over 
    r(1/\epsilon w) }
  {r(1/\epsilon w)\over r(1/\epsilon) }
  \hfill\cr\hfill
  {}+ \lim_{\epsilon\to 0} 
  {F^\leftarrow(1-\epsilon wy)-F^\leftarrow (1-\epsilon y) \over r(1/\epsilon)}
  \, .
  \qquad
  \hbox{(A.1)}\cr
  }
$$
It follows that $\lim_{\epsilon\to 0} r(1/\epsilon w)/r(1/\epsilon)$
exists for almost all $w$. Hence $r$ is regularly varying and there exists
a real number $\rho$ such that $\lim_{\epsilon\to 0} r(1/\epsilon w)/
r(1/\epsilon)= w^\rho$. Then (A.1) yields the functional equation
$$
  \tilde h(uw)=\tilde h(u)w^\rho + \tilde h (w) \, .
$$

If $\rho$ vanishes, this means $\tilde h(uw)=\tilde h(u)+\tilde h(w)$.
Since $\tilde h$ is monotone, this forces it to be proportional to the 
logarithm function.

If $\rho$ does not vanish, then, permuting $u$ and $w$, we obtain
$$
  \tilde h(wu)=\tilde h(w)u^\rho + \tilde h (u) \, .
$$
Hence, equating the expressions obtained for $h(uw)$ and $h(wu)$, we have
$$
  \tilde h(u)(w^\rho-1)=\tilde h(w) (u^\rho-1)\, .
$$
This implies that the function $\tilde h(u)/(u^\rho-1)$ is almost everywhere
constant. Thus, there exists a constant $c$ such that $\tilde h=ck_\rho$
almost everywhere. Again, 
since $\tilde h$ is monotone, this almost everywhere equality holds in fact
everywhere.

In any case, regardless whether $\rho$ vanishes or not, we obtain that
$\tilde h=ck_\rho$ for some constant $c$. This means,
setting $c_1=h(x)$ and $c_2=c\bigl( h(v)-h(y)\bigr)$,
$$
  h(u)=c_1+c_2k_\rho(u)\,. 
$$
The function $h$ is then continous on the positive half line. Therefore,
Lemma 2.3.ii shows that $h_{1/\epsilon}$ converges to $h$ everywhere as
$\epsilon$ tends to $0$. In particular,
$$
  \lim_{\epsilon\to 0} 
  { h_{1/\epsilon}(u)-h_{1/\epsilon}(1)\over
    h_{1/\epsilon}(v)-h_{1/\epsilon}(1) }
  = {u^\rho-1\over v^\rho-1} \, .
$$

The same argument applies for what is needed in the proof of the 
sufficiency part of Theorem 3.1, namely that if
$$
  \lim_{\epsilon\to 0} 
  { F^\leftarrow(1-\epsilon u)-F^\leftarrow (1-\epsilon x)
    \over F^\leftarrow(1-\epsilon v)-F^\leftarrow (1-\epsilon y) }
$$
exists for almost all $u,v,x,y$, then it exists for all $u,v,x$ and $y$.
This comes from the fact that the functions
$$
  u\mapsto 
  { F^\leftarrow (1-\epsilon u)-F^\leftarrow(1-\epsilon x)
    \over
    F^\leftarrow (1-\epsilon v)-F^\leftarrow(1-\epsilon y) }
$$
are monotone and that one can take $v,x$ and $y$ such that these functions
converge for almost all $u$ as $\epsilon$ tends to $0$; hence, the
convergence occures for all $u$ and locally uniformly; permuting the
variables $u,v,x$ and $y$, the convergence also occurs locally uniformly
with respect to all $u,v,x$ and $y$.\hfill\qed

}

\bigskip

\noindent{\bf Acknowledgements.} Years ago, Anne-Laure Foug\`eres kept asking
me questions about the linear normalization in asymptotic extreme value
theory. It is a pleasure to acknowledge that her questions are at the root
of this note, and that my students at the Universit\'e de Cergy-Pontoise
during the spring 2008 term for prompting me to write the proof in 
this paper. This
note also benefited from comments, remarks and suggestions from Bill 
McCormick, precise, constructive, numerous, needed and welcome as always.

\bigskip


\noindent{\bf References}
\medskip

{\leftskip=\parindent \parindent=-\parindent
 \par

P.\ Billingsly (1968). {\sl Weak Convergence of Probability Measures}, 
Willey.

N.H.\ Bingham, C.M.\ Goldie, J.L.\ Teugels (1989). {\sl Regular Variation},
Cambridge.

L.\ de Haan (1970). {\sl On Regular Variation and its Application to the Weak
Convergence of Sample Extremes}, Mathematical Centre Tracts 32, Mathematisch 
Centrum Amsterdam

W.\ Feller (1970). {\sl An Introduction to Probability Theory and its 
Applications}, Wiley.

E.\ Hlawka, J.\ Schoi\ss engeier, R.\ Taschner (1986). {\sl Geometric and
Analytic Number Theory}, Springer.

M.R.\ Leadbetter, G.\ Lindgren, H.\ Rootz\'en (1983). {\sl Extreme and Related
Properties of Random Sequences and Processes}, Springer.

S.I.\ Resnick (1987). {\sl Extreme Values, Regular Variation, and Point
Process}, Springer.

W.\ Rudin (1976). {\sl Principles of Mathematical Analysis}, McGraw-Hill.

}


\bigskip

\setbox1=\vbox{\halign{#\hfil\cr
  Ph.\ Barbe            \cr
  90 rue de Vaugirard   \cr
  75006 PARIS           \cr
  FRANCE                \cr
  philippe.barbe@math.cnrs.fr\cr}}
\box1

\bye